\documentclass[12pt,a4paper]{article}
\usepackage{localeng}
\usepackage{hyperref}
\usepackage[margin=20mm]{geometry}
\usepackage{enotez}
\DeclareInstance{enotez-list}{custom}{paragraph}
{
 heading = \section*{#1} ,
 notes-sep = \baselineskip ,
 format = \normalfont ,
 number = \enmark{#1}
}
\usepackage{comment}

\newcommand{\cnd}{\mskip 1mu|\mskip 1mu}

\begin{document}

\title{Kolmogorov's Last Discovery?\\ (Kolmogorov and algorithmic statistics)}
\author{Alexey Semenov\thanks{Moscow State Lomonosov University, Moscow, Russia}, Alexander Shen\thanks{LIRMM, Univ Montpellier, CNRS, Montpellier, France. Supported by 
ANR-21-CE48-0023 FLITTLA grant.}, Nikolay Vereshchagin\thanks{Moscow State Lomonosov University, Moscow, Russia, National Research University The Higher School of Economics (HSE), Moscow, Russia. The article was supported by the Fundamental Research program of HSE.}}
\date{}
\maketitle

The last theme of Kolmogorov's mathematics research was \emph{algorithmic theory of information}, now often called \emph{Kolmogorov complexity theory}. Kolmogorov played a crucial role in its creation, though he was not the first who suggested measuring the amount of information in a finite object using theory of algorithms, as the length of the minimal description. Ray Solomonoff did this several years earlier, in the beginning of 1960s.\endnote{In Solomonoff's report written in 1960~\cite{1960solomonoff} he wrote:
\begin{quote}
Consider a very long sequence of symbols --- e.g., a passage of English text, or a long mathematical derivation. We shall consider such a sequence of symbols to be ``simple'' and have a high a priori probability, if there exists a very brief description of this sequence --- using, of course, some sort of stipulated description method. More exactly, if we use only the symbols $0$ and $1$ to express our description, we will assign the probability $2^{-N}$ to a sequence of symbols, if its shortest possible binary description contains $N$ digits.
\end{quote}

The notion of ``binary description'' is later explained as follows:
\begin{quote}
Suppose that we have a general purpose digital computer $M_1$ $\langle\ldots\rangle$ 

Any finite string of $0$'s and $1$'s is an acceptable input to $M_1$. The output of $M_1$ (when it has an output) will be a (usually different) string of symbols, usually in an alphabet other than the binary. If the input string $S$ to machine $M_1$ gives output string $T$, we shall write
\[
M_1(S)=T.
\]
Under these conditions, we shall say that ``$S$ is a description of $T$ with respect to machine $M_1$''.  If $S$ is the shortest such description of $T$, and $S$ contains $N$ digits, then we will assign to the string $T$ the a priori probability $2^{-N}$.
\end{quote}
In the same report Solomonoff suggested to consider a universal machine that can simulate every machine if the appropriate prefix is added to the input. Such a machine provides the shortest (up to an additive $O(1)$ term) description. But Solomonoff' s main goal was to introduce a natural notion of \emph{a priori} probability distribution on finite objects, and he immediately bumped into a technical problem. He suggested assigning probability $2^{-N}$ to an object whose shortest description contains $N$ bits. But if we add the values of $2^{-N}$ for all shortest descriptions, we may get a sum greater than $1$. Solomonoff suggested some (ad hoc) way to deal with this problem, but mentioned that they are not very convincing. In his 1964 paper~\cite{1964-1solomonoff} he writes: 

\begin{quote}
The author feels that the proposed systems are consistent and meaningful, but at the present time, this feeling is supported only by heuristic reasoning and several nonrigorous demonstrations. [p.~4]
\end{quote}

And later (in the same paper):
\begin{quote}
The author feels that Eq. (1) [some way to overcome the problem of the total probability greater than $1$] is likely to be correct or almost correct, but that the methods of working with the problems of Sections 4.1 to 4.3 are \emph{more} likely to be correct than Eq.~(1). If Eq.~(1) is found to be meaningless, inconsistent, or somehow gives results that are intuitively unreasonable, then Eq.~(1) should be modified in ways that do not destroy the validity of the methods used in Sections 4.1 to 4.3 [p.~10].
\end{quote}

One can easily imagine the mathematicians' reactions to passages like this. But in fact Solomonoff's idea of a priori probability that appears when we run a random program, and its connections with the notion of complexity (the length of a shortest description) turned out to a be a key point in algorithmic information theory both from philosophical and technical viewpoints. Here is how Solomonoff describes the differences between these two notions:

\begin{quote}
Another objection to the method outlined is that Equation (1) uses only the ``minimal binary description'' of the sequences it analyzes. It would seem that if there are several different methods of describing a sequence, each of these methods should be given \emph{some} weight in determining the probability of that sequence. [p.~16]
\end{quote}
}~\cite{1960solomonoff,1964-1solomonoff,1964-2solomonoff}. However, Solomonoff's work remained almost unnoticed before Kolmogorov's publications (the first was in 1965) where he referred to Solomonoff's work. Probably this happened not only because of Kolmogorov's fame but also because Kolmogorov papers were well written with clear definitions and proofs (and informal motivation, too). 

There were two main publications of Kolmogorov related to algorithmic information theory. The first appeared in 1965~\cite{1965kolmogorov}, and the second one~\cite{1969kolmogorov} appeared in 1968 (the English version) --- 1969 (the Russian version); its results were presented in September 1967 at the International Symposium on Information Theory (San Remo, Italy) and in October 1967 at Moscow Mathematical Society meeting\endnote{%
Here is the translation of the abstract of this talk~\cite{1968umn}.
\begin{quote}
\textbf{October 31, 1967}

\textbf{1}. \textit{A.N.~Kolmogorov} Several results on algorithmic entropy and algorithmic amount of information.

The algorithmic approach to the foundations of information theory and probability theory was obstructed by some unsolved questions that appeared from the very beginning. Now the situation has changed. In particular, it was shown that the decomposition formula $H(x,y)\sim J(x)+H(y\cnd x)$ as well as the formula $J(x\cnd y)\sim J(y\cnd x)$ are true in the algorithmic framework only up to $O[\log H(x,y)]$ additive terms (Levin, Kolmogorov).

The statements about fundamental difference between the algorithmic definition of Bernoulli sequences (the simplest case of Kollektiv) and Mises--Church definition are now proven in the following form: there are sequences  $x=(x_1,x_2,\ldots)$ that satisfy the Mises--Church definition with density $1/2$ of ones whose initial segments have entropy (complexity) $H(x^n)=H(x_1,x_2,\ldots,x_n)=O(\log n)$ (Kolmogorov).  The talk is accessible to the audience familiar with the informal notion of a computable function.
\end{quote}
Here Kolmogorov uses the name ``entropy'' instead of complexity and denotes it $H$ instead of $K$. The letter $J$ stands for mutual information.
} 

In addition to these two papers, Kolmogorov's survey paper~\cite{1970kolmogorov}  based on his Nice International Congress of Mathematicians was published (with big delay, only in 1983). Another survey paper based on the Kolmogorov--Uspensky's talk at the congress of the Bernoulli Society (1987, ~cite{1987kolmogorov}) was prepared mostly by Uspensky alone (as noted by Uspensky) --- at that time Kolmogorov was very ill. So Kolmogorov's ideas that did not appear as proved (and published) theorems can be reconstructed only partially based on work of his students and collaborators, short abstracts of his talks and the recollections of people who were present at these talks.

In this survey we try to reconstruct the development of Kolmogorov's ideas related to algorithmic statistics (resource-bounded complexity, structure function and stochastic objects).

\section{Resource-bounded complexity}

The 1965 paper~\cite{1965kolmogorov} the definition of a complexity of a finite object is given (as the minimal length of program generating this object), and it is shown that there exist optimal programming languages that make the complexity minimal up to $O(1)$ additive term. In addition, this paper contained some remarks in its closing section (Final remarks). One of them was the following:

\begin{quote}
The approach of Section 3 has one important deficiency: it does not take into account how ``hard'' is to transform the program $p$ and object $x$ into object $y$. Giving appropriate definitions, one can formulate and prove rigorously some mathematical results that can be  reasonably interpreted as indications that there are some cases when an object has a very short program (and therefore has very small complexity $K(x)$, but can be reconstructed from short programs only by unrealistically long computations. I plan to study elsewhere the dependence of the complexity (the minimal program length) $K^t(x)$ on the allowed difficulty $t$ of transforming the program into~$x$. Then the complexity $K(x)$ introduced in Section 3, reappears as the minimum of $K^t(x)$ for unlimited values of $t$.
\end{quote}

This paragraph essentially contains the definition of the resource-bounded complexity. For example, for the time-bounded case we say that $K^t(x)$ is the minimal length of a program that produces $x$ in at most $t$ steps. When $t$ increases, the value of $K^t(x)$ decreases and (starting from some moment) becomes equal to the usual (unbounded) complexity $K(x)$.

A cautious remark about what can be done after giving proper definitions and how one could interpret these results, can be formalized in the following way. Let $t(n)$ be an arbitrary computable total function with natural arguments and values. (The interesting case is when $t(n)$ grows fast as $n$ grows). Then for every $n$ one can find a string $x$ of length $n$ such that  $K^{t(n)}(x_n)\ge n$ while $K(x_n)\le O(\log n)$. Indeed, we may consider the lexicographically first string $x$ of length $n$ such that $K^{t(n)}(x)\ge n$. Such a string exists (even the string of length $n$ with $K(x)\ge n$ exists). This string can be computed given $n$ by simulating all the programs of length at most $n$ for $t(n)$ step. This is a lengthy process, but for resource-unbounded complexity only computability matters. Therefore, the (unbounded) complexity of $n$ is $O(\log n)$.

For every string $x$, the complexity $K(x)$ is a numerical quantity. The resource-bounded complexity of some fixed $x$ is no more a number, it is a (decreasing) function $t\mapsto K^t(x)$. Kolmogorov (see the quote above) promised to study this function elsewhere, but he did not publish other papers about time-bounded complexity.\endnote{%
In 1965 Kolmogorov gave a popular talk in the Institute of Philosophy (Academy of Sciences, USSR). A transcript (almost unedited) of this talk was found among Kolmogorov's papers after his death and published in~\cite{1965kolmogorovtalk}. In this talk he said:
\begin{quote}
One can show (this theme will be developed in the paper for \emph{Uspekhi matematicheskikh nauk} that I am preparing) [such a paper was never published; the publishers of the transcript wrongly assumed that Kolmogorov had in mind the paper~\cite{1970kolmogorov}, but the latter was a survey prepared for Nice IUM congress] that there are cases where, for example, a problem with a simple statement has a solution that can be symbolically expressed by a rather short formula. But at the same time it is known that if we want to avoid an unreasonably large amount of computational work, the solution has only much longer representations. Mathematicians are now able to prove results of this type, and these theorems should in fact replace G\"odel's [incompleteness] theorem which is often mentioned.
\end{quote}
}
It is possible that this happened because of two types of problems related to resource-bounded complexity.
\begin{itemize}
\item The definition of bounded-resource complexity is machine-dependent, since the resource usage changes when we change the computation model, and there is no ``right'' computational model.
\item Trying to prove something about the time-bounded complexity (which seems to be the most natural and interesting approach), we see almost immediately that the answer may depend on the P=NP question (or related ones).
\end{itemize}

Both problems disappear if we measure the resources in a very rough way --- so rough that the differences between models do not matter anymore. Kolmogorov noted in his talk at another meeting of the Moscow Mathematical society in 1972 (as the published abstract\endnote{%
Here is the full translation of this abstract~\cite{1972umn}:
\begin{quote}
\textbf{November 23, 1971}

\textbf{1}. \emph{A.N.~Kolmogorov}, The description complexity and the construction complexity of mathematical objects.

$1^\circ$. When designing the computer programs we have to estimate (a)~the complexity of a program; (b)~the space used by the program; (c)~the duration of the computation. In this talk we consider a more theoretical analysis of this type of question.

$2^\circ$. In 1964--1965 it was noted that the minimal length $K(x)$ of a program (in binary) that constructs some object $x$ can be defined in an invariant way up to $O(1)$ additive terms (Solomonoff, Kolmogorov).  This allows to consider the notion of \emph{description complexity} of a constructive mathematical object as a starting point for a new approach to the foundation of information theory (Kolmogorov, Levin) and probability theory (Kolmogorov, Martin-L\"of, Schnorr, Levin).

$3^\circ$. The situation becomes more difficult if we want to have a machine-independent measure of the space/time required by a computation. But some results can be obtained in the framework of machine-independent computation complexity theory (Blum, 1967). Let $\Pi(p)$ be some characterization of ``complexity'' of the construction of an object $x=A(p)$ according to a program $p$, and let $\Lambda(p)$ be the length of the program $p$. Then we define the ``$n$-complexity'' of an object $x$ as $K^n\Pi(x)=\inf\,\{\Lambda(p): x=A(p), \Pi(p)=n\}$ [a typo: should be $\Pi(p)\le n$ instead] (if the conditions cannot be fulfilled, the infimum is infinite).

$4^\circ$. Barzdins' theorem about the complexity $K(M_a)$ of initial segments of a indicator sequence of a [computably] enumerable set of natural numbers (1986) and the results of Barzdins, Kanovich and Petri about the corresponding complexities $K^n\Pi(M_a)$ are of general mathematical interest and provide some new view on the role of exceeding the limits of the formalized part of mathematics that was used before. The survey of these results was made without tedious technical details.
\end{quote}
}
says) that

\begin{quote}
The situation becomes more difficult if we want to have a machine-independent measure of the space/time required by a computation. But some results can be obtained in the framework of machine-independent computation complexity theory (Blum, 1967). Let $\Pi(p)$ be some characterization of ``complexity'' of the construction of an object $x=A(p)$ according to a program $p$, and let $\Lambda(p)$ be the length of the program $p$. Then we define the ``$n$-complexity'' of an object $x$ as $K^n\Pi(x)=\inf\,\{\Lambda(p): x=A(p), \Pi(p)=n\}$ [a typo: should be $\Pi(p)\le n$ instead] (if the conditions cannot be fulfilled, the infimum is infinite).
\end{quote}

No further details were provided.

As to the second obstacle, Kolmogorov understood from the very beginning the importance of P=NP question (and encouraged Levin to publish his paper about NP-completeness), and definitely understood the obstacles related to this question in the bounded-resource complexity theory.

Longpr\'e and Mocas write in their 1993 paper~\cite{1993longpre} that ``As Kolmogorov stated himself (Levin, Private communication), the problem of whether symmetry holds in time bounded environment has interesting connections to complexity theory''. In this paper they show that the symmetry of information (formula for the complexity of a pair) for time-bounded complexity does not hold if one-way functions exist. They also present (following Longpr\'e's PhD thesis~\cite{1986longpre}) the proof of the symmetry of information for the space-bounded complexity.

Anyway, in 1980s several versions of time-bounded complexity were known and used --- for example, Sipser in~\cite{1983sipser} used time-bounded complexity to show that the complexity class BPP belongs to polynomial hierarchy (this result was soon improved to $\Sigma_2 \cap \Pi_2$ by G\'acs, see the same Sipser's paper).

Earlier, in the 1970s, Bennett pointed out that small complexity can coexist with large resource-bounded complexity. Chaitin in 1977 paper~\cite{1977chaitin} quotes Bennett:

\begin{quote}
C.H.~Bennett [55] has suggested another approach [to give a quantitative structural characterization\ldots of degree of organization] based on the notion of ``logical depth''. A structure is deep ``if it is superficially random but subtly redundant, in other words, if almost all its algorithmic probability is contributed by slow-running programs. A string's logical depth should reflect the amount of computational work required to expose its buried redundancy''. It is Bennett's thesis that ``a priori the most probable explanation of `organized information' such as a sequence of bases in a naturally occurring DNA molecule is that it is the product of an extremely long evolutionary process''. $\langle\ldots\rangle$
We are grateful to C.H.~Bennett for permitting us to present his notion of logical depth in print for the first time\ldots
\end{quote}

Here [55] is a reference to an undated unpublished manuscript of Bennett ``On the Thermodynamics of Computation''.

In 1979 Adleman~\cite{1979adleman} mentions a specific example when the (informal) description complexity decreases when the allowed computational resources increase. He notes that (a)~the number $2^{136}+1$ has a short description and is easy to compute (starting from this description); (b)~a random number of similar size does not have short descriptions at all; (c) the number that is \emph{the maximal prime factor of $2^{136}+1$} has a short description (just given), but a lot of efforts are required to reconstruct this number starting from this description (factorization is difficult).

In 1988 Bennett himself published a paper~\cite{1988bennett} in a collection devoted to the 50th anniversary of universal Turing machines when he says\endnote{This quote is taken from the Abstract; in the paper itself Bennetts gives more details:
\begin{quote}
A string is said to be \emph{compressible} by $k$ bits if its minimal program is $\ge k$ bits shorter than the string itself. A simple counting argument shows that at most a fraction $2^{-k}$ of strings of length $\le n$ bits can have this property. This fact justifies calling strings that are incompressible, or nearly so, \emph{algorithmically random}. Like the majority of strings one might generate by coin tossing, such strings lack internal redundancy that could be exploited to encode them concisely $\langle\ldots\rangle$ 

The relation between universal computer programs and their outputs has long been regarded [40] as a formal analog of the relation between theory and observation in science, with the minimal-sized program representing the most economical, and therefore a priori most plausible, explanation of its output. This analogy draws its authority from the ability of universal computers to execute all formal deductive processes and their presumed ability to simulate all processes of physical causation. Accepting this analogy, one is then led to accept the execution of the minimal program as representing its output’s most plausible causal history, and a logically ``deep'' or complex object would then be one whose most plausible origin, via an effective process, entails a lengthy computation. Just as the plausibility a scientific theory depends on the economy of its assumptions, not on the length of the deductive path connecting them with observed phenomena, so a slow execution time is not evidence against the plausibility of a program; rather, if there are no comparably concise programs to compute the same output quickly, it is evidence of the nontriviality of that output.

A more careful definition of depth should not depend only on the minimal program, but should take fair account of all programs that compute the given output, for example giving two $k+1$ bit programs the same weight as one $k$-bit program.
\end{quote}
Here [40] is Solomonoff's paper~\cite{1964-1solomonoff}.
 }%
that
\begin{quote}
some mathematical and natural objects (a random sequence, a sequence of zeros, a perfect crystal, a gas) are intuitively trivial, while others (e.g. the human body, the digits of $\pi$) contain internal evidence of a nontrivial causal history. We formalize this distinction by defining an object’s ``logical depth'' as the time required by a standard universal Turing machine to generate it from an input that is algorithmically random\ldots
\end{quote}

Looking for a suitable formalization, Bennett starts with the following definition of \emph{depth}:
\begin{quote}
\textbf{Tentative Definition 0.1}: A string’s depth might be defined as the execution time of its minimal program.

The difficulty with this definition arises in cases where the minimal program is only a few bits smaller than some much faster program, such as a print program, to compute the same output $x$. In this case, slight changes in $x$ may induce arbitrarily large changes in the run time of the minimal program, by changing which of the two competing programs is minimal. Analogous instability manifests itself in translating programs from one universal machine to another.
\end{quote}

A more robust approach, says Bennett, is to consider the dependence between the allowed computation time and the length of program (as Kolmogorov suggested): 
\begin{quote}
\textbf{Tentative Definition 0.2}: A string’s depth at significance level $s$ [might] be defined as the time required to compute the string by a program no more than $s$ bits larger than the minimal program.
\end{quote}

Then Bennett considers further modification of this definitions, but we will not consider this modification here (as well as his definition of depth for infinite sequences).

Similar definitions appear in many subsequent papers. For example, in~\cite{2006AFMV} the \emph{computational depth} is defined as follows (time is considered as an independent variable; the Kolmogorov complexity is denoted by letter $C$, and $C^t$ is complexity with time bound~$t$):
\begin{quote}
\textbf{Definition 4.1} Let $t$ be a time bound. The time-$t$ depth of $x$ is 
\[
D^t(x) = C^t(x) - C(x).
\]
\end{quote}

Let us return to the question of machine-dependence in the definition of resource-bounded complexity. This dependence becomes negligible if we consider ``astronomically large'' time bounds (or space bounds, the difference does not matter anymore) up to arbitrary computable transformations. Technically this can be done by considering the function $B(n)$ defined as the maximal number that has complexity at most $n$ (as suggested in~\cite[Section 2.2]{1983gacs}). Another equivalent (up to $O(1)$ change in the argument) definition of this function: the maximal computation time for terminating programs of size at most $n$. (If we consider the number of states of Turing machines instead of program size, we get ``busy beaver numbers''.)

Now, for a given finite object $x$, we consider a function 
\[
k\mapsto K^{B(k)} (x)
\] 
It decreases as $k$ increases and reaches its final value around $K(x)$ (or earlier), since the time required for the termination of the minimal program for $x$ is bounded by $B(K(x))$. It turned out that the same (with logarithmic precision) function appears in two other contexts considered by Kolmogorov: for two-part descriptions and for stochastic objects. These notions are considered in the two following sections.

\section{Two-part descriptions and Kolmogorov's structure function}
\label{sec:structure}

Consider an asymmetric coin when one side is (for example) twice more probable than the other one. Toss this coin $n$ times and write the result as an $n$-bit sequence. When $n$ is large enough, such a sequence will (with overwhelming probability) contain approximately $1/3$ of zeros and $2/3$ of ones. (We assume that $1$ means the more probable side). The complexity of this sequence will be significantly less than $n$; this makes the situation different from the case of a symmetric coin. Our sequence can be described using two pieces of information:
\begin{itemize}
\item first we specify the number of zeros and ones in the sequence;
\item then we specify the ordinal number of our sequence in the lexicographically ordered list of all sequences with the same number of zeros and ones.
\end{itemize}
The first part requires $O(\log n)$ bits, and the second part requires about $\log \binom{n}{n/3}$ bit. This number is close to $nH(1/3)$, where 
\[
H(p) = p \log\frac{1}{p}+ (1-p)\log\frac{1}{1-p}
\]
is the Shannon function. So the first part is negligible compared to the second part. With high probability (with respect to the Bernoulli distribution we consider) this two-part description is close to the optimal, and the Kolmogorov complexity of the sequence is close to $nH(1/3)$.

These two-part descriptions appear in more general situations. Consider some finite object $x$ (the outcome of some experiment) that is contained in some finite set $A$ of finite objects. Then we may consider the following \emph{two-part description} of $x$:
\begin{itemize}
\item the first part is a description of $A$ and uses $K(A)$ (if we use the optimal description);
\item the second part is an ordinal number of $x$ in $A$ and uses $\log_2\#A$ bits.
\end{itemize}
Here we assume that the finite set $A$ is encoded in some natural way, say, as a binary string that allows us to reconstruct the list of all elements of $A$, and $K(A)$ denotes the complexity of this binary string (the choice of a computable encoding does not matter with $O(1)$-precision). We also fix some natural ordering of $A$ for the second part (say, the lexicographical ordering if elements of $A$ are binary strings).

How close is this two-part description to being optimal? In other words, how much the sum $K(A)+\log_2\#A$ exceeds $K(x)$? Note that this sum cannot be smaller than $K(x)$ since the optimal description should be optimal. Of course, the answer depends on $A$ and $x$. In our example (asymmetric coin) we get a description that is (with high probability) close to the optimal one. This fact can be interpreted as follows: our set $A$ is a good model for the results of the coin-tossing experiment. 

We could consider the set $A$ of all $n$-bit sequences instead; then we get a two-part description of size $O(\log n)+n$, since now the second part would require $n$ bits. This description is not optimal, and, indeed, our ``model'' does not take into account some specific property of $x$ (that says that frequency of ones is close to $2/3$).

On the other hand, for every $x$ we may consider a singleton $A=\{x\}$; its complexity is the same as for $x$, but $\log_2\#A=0$, so we get an optimal two-part description. Being optimal, this description, informally speaking, is ``too specific'' --- we now consider some random peculiarities of $x$ as part of the model. It would be better to avoid this, in other words, to choose (among all optimal two-part description) a description with small first part (and large second part).

This approach was description in Kolmogorov's talk on the April 16, 1974 meeting of Moscow Mathematical Society (according to the published abstract~\cite{1974umn}):
 \begin{quote}
 \textbf{1}. \emph{A.N.~Kolmogorov}. The complexity of algorithms and the objective definition of randomness.
 
 For every constructive object $x$ we consider a function $\Phi_x(k)$ whose argument is a natural number $k$ and whose value is the logarithm of a minimal size of a finite set of complexity at most $k$ that contains $x$. If $x$ has a simple description, the function $\Phi_x(k)$ equals $1$ [typo: should be $0$] even for small values of $k$. If $x$ does not have a simple description, the object $x$ is ``random in the negative sense''. But $x$ is ``positively probabilistically random'' only if the function $\Phi$ reaches some value $\Phi_0$ for some relatively small $k=k_0$ and then decreases approximately as $\Phi(k)=\Phi_0-(k-k_0)$.

 \end{quote}
 
 It seems that Kolmogorov mentioned the same approach in his talk at the international conference on information theory in Tallinn (1973). A picture of Kolmogorov giving a talk was made by Fine~\cite{1973fine};
\begin{center}
\includegraphics[width=0.9\textwidth]{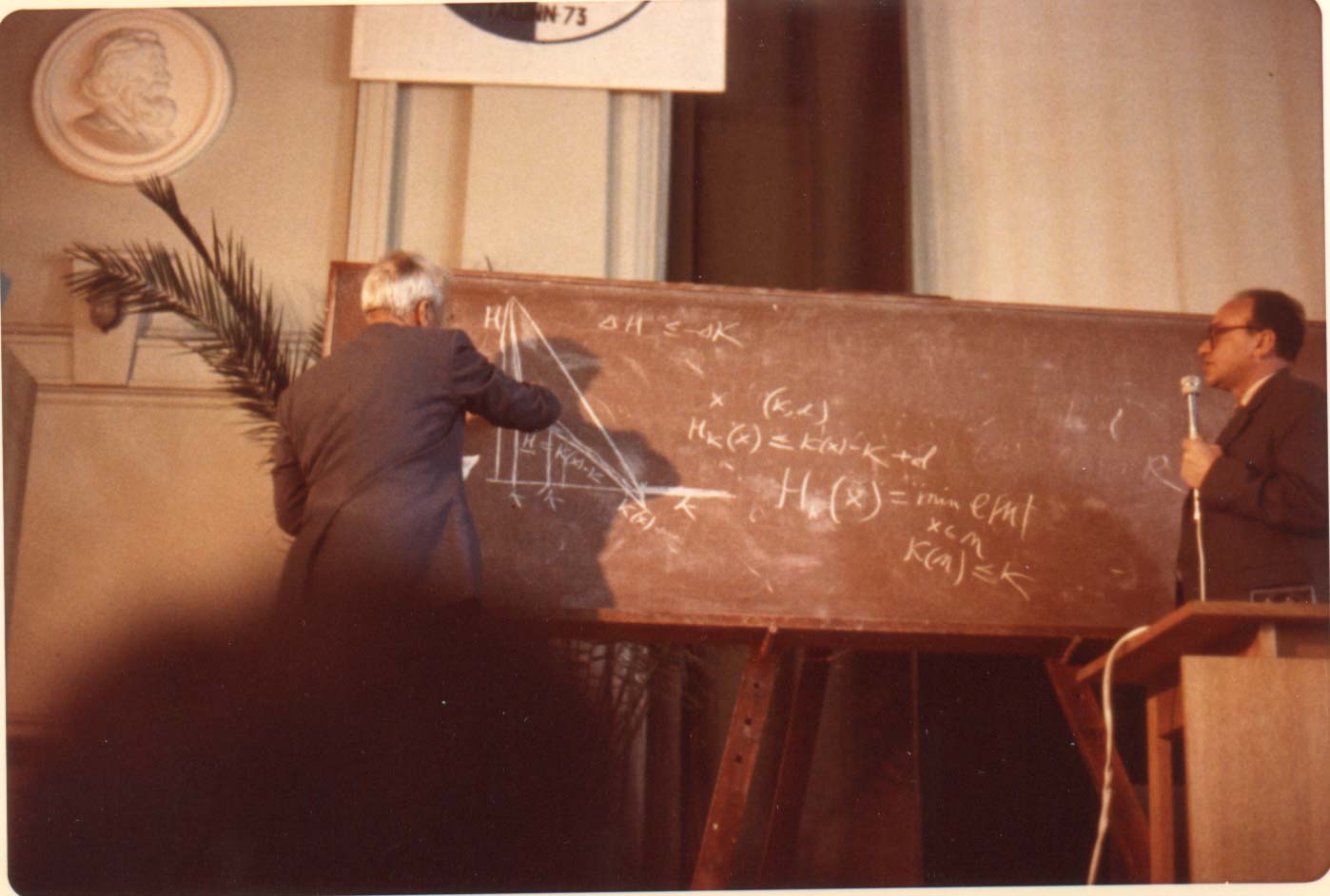}
\end{center}
we see the same definition on the board (though $\Phi$ is replaced by $H$).

Unfortunately, we do not have first-hand information about the contents of Kolmogorov's talk except for this photo\endnote{There is an audio recording of \emph{some} talk by Kolmogorov given at some international conference~\cite{1973dynkin}. It may be the first part of Tallinn's talk, but the talk starts slowly and the recorded part does not say anything about the $H$ function. So it may be some other talk given elsewhere. (It cannot be the last talk by Kolmogorov given in Tbilisi~\cite{1983kolmogorovconf}, since at the latter conference the translation interrupted Kolmogorov's speech quite often.)
}
 
Here is the description of what Kolmogorov said in Tallinn according to~\cite[p.~32]{1985cover}:\endnote{Another account of what Kolmogorov said is given in~\cite[p.~175--176]{1991CT}:
\begin{quote}
We begin with a definition of the smallest set containing  $x^n$ that is describable in no more than $k$ bits.

\textit{\textbf{Definition}}. The \emph{Kolmogorov structure function} $K_k(x^n\cnd n)$ of a binary string $x\in \{0,1\}^n$ is defined as
\[
K_k(x^n\cnd n)=\min_{p: l(p)\le k;\  \mathcal{U}(p,n)=S,\  x^n\in S\subseteq\{0,1\}^n} \log |S| \eqno (7.104)
\]
The set $S$ is the smallest set which can be described with no more than $k$ bits and which includes $x^n$. By $\mathcal{U}(p,n)=S$, we mean that running the program $p$ with data $n$ on the universal computer $\mathcal{U}$ will print out the indicator function of the set $S$.

\textit{\textbf{Definition}}.  For a given small constant $c$, let $k^*$ be the least $k$ such that
\[
K_k(x^n\cnd n) + k \le K(x^n\cnd n)+c\eqno(7.105)
\]
Let $S^{**}$ be the corresponding set and let $p^{**}$ be the program that prints out the indicator function of $S^{**}$. Then we shall say that $p^{**}$ is a \emph{Kolmogorov minimal sufficient statistic} for $x^n$.

\begin{center}
$\langle\ldots\rangle$
\end{center}
The concept of the Kolmogorov structure function was defined by Kolmogorov at a talk at the Tallin[n] conference in 1973, but these results were not published. 
\end{quote}
}
 
\begin{quote}
Consider the function $H_k\colon\{0,1\}^n\to N$, 
\[
H_k(x)= \min_{p\colon l(p)\le k} \log |S|,
\]
where the minimum is taken over all subsets $S\subseteq \{0,1\}^n$, such that $x\in S$, $U(p)=S$, $l(p)\le k$. This definition was introduced by Kolmogorov in a talk at the Information Theory Symposium, Tallin[n], Estonia, in 1974.\footnote{Most probably an error: the published volume of abstracts gives the date as 1973, June 18--23.} Thus $H_k(x)$ is the log of the size of the smallest set containing $x$ over all sets specifiable by a program of $k$ or fewer bits. Of special interest is the value
\[
k^*(x) = \min \{k\colon H_k(x)+k = K(x)\}.
\]
Note that $\log |S|$ is the maximal number of bits necessary to describe an arbitrary element $x\in S$. Thus a program for $x$ could be written in two stages: ``Use $p$ to print the indicator function for $S$; the desired sequence is the $i$th sequence in a lexicographic ordering of the elements of this set.'' This program has length $l(p)+\log|S|$, and $k^*(x)$ is the length of the shortest program $p$ for which this $2$-stage description is as short as the best $1$-stage description $p^*$. We observe that $x$ must be maximally random with respect to $S$ --- otherwise the $2$-stage description could be improved, contradicting the minimality of $K(x)$. Thus $k^*(x)$ and its associated program $p$ constitute a minimal sufficient description for $x$.
\end{quote} 

However, there are other accounts of what Kolmogorov said in Tallinn (see in the next section about stochastic objects).
\smallskip

The two-part descriptions were considered (in a slightly different form) by other authors. In 1988 Koppel~\cite{1988koppel} describes this approach as follows:
 
\begin{quote}
3. What is the sophistication of the string, that is, what is the minimal amount of planning which must have gone into the generation of the string? More picturesquely, if the string is being broadcast by some unknown source, what is the minimal amount of intelligence we must attribute to that source? [p.~435]

\ldots both simple strings and random strings are not sophisticated, in the sense explained in Question 3. $\langle\ldots\rangle$

The relationship between sophistication and complexity can be made more precise in the following way. 

The minimal description of a strings consists of two parts. One part is a description of the string's structure, and the other part specifies the string from among the class of strings sharing this structure (Cover 1985). The sophistication of a string in the size of that part of the description which describes the string's structure. Thus, for example, the description of the structure of a random string is empty and thus, though its complexity is high, its sophistication is low.

$\langle \ldots \rangle$

Having defined complexity in terms of a two-part description, it is easy to sort out the ``sophistication'' from the complexity.

\textbf{Definition}. \emph{The $c$-sophistication of $S$}, 
\begin{multline*}
\textit{SOPH}_c(S) = \min\{|P|\ \mid\ \exists D\, \text{$(P,D)$ \emph{is a description of $S$}}\\ 
\text{\emph{and $|P|+|D|\le H(S)+c$}}\}.
\end{multline*}
 We call a description $(P,D)$ of $S$, \emph{$c$-minimal} if $|P|+|D|\le H(S)+c$. We call a program $P$ a \emph{$c$-minimal program} for $S$ if $P$ is the shortest program such that for some $D$, $(P,D)$ is a $c$-minimal description of $S$. Thus $\textit{SOPH}_c(S)$ is the length of a $c$-minimal program for~$S$.

The $c$-minimal program for $S$ is that part of the description of $S$ which compresses $S$ --- it represents the structure of $S$. The range of this program constitutes the class of strings which share the structure of $S$.
\end{quote}
 Here ``(Cover 1985)'' is the paper~\cite{1985cover} mentioned above. The pair $(P,D)$ consists of a program $P$ and some input $D$, and we require that $P$ is a total program (that terminates on all inputs). Then we may consider the finite set $A$ of all $P(D')$ for all $D'$ that have the same length as $D$ (``the range of this program''). The program $P$ (together with the binary representation of the length of $D$, but the latter is quite short and we ignore it) determines the set $A$,  and the length of $D$ is the binary logarithm of the size of $A$, so we get almost the same definition as before. There are technical differences, though (except for the length of $D$ that was already mentioned): Koppel assumes that every description of some binary string $S$ is at the same time a description of all prefixes of $S$. Also he considers the curve in different coordinates: he starts from some $c$ which is an allowed overhead when we switch from the optimal description of $S$ (whose length is denoted by $H(S)$) to two-part descriptions of the same $S$ and looks how large should be the length of the first part of the two-part description with this (or smaller) overhead.  (See also other Koppel's publications that appeared around the same time~\cite{1987koppel,1991koppelatlan}.)
 
One can describe this approach in very general (almost philosophical, cf. the survey~\cite{1992li-vitanyi}) terms. Science constructs models that are used to describe the  results of observations and experiments. The model should be significantly simpler that the data it describes: if the model has more parameters than there are experimental points, it is hardly interesting (``overfitting''). One can remember Occam's razor and Mach's ``economy of thought'' principle.\endnote{Readers who lived during Soviet times also remember Lenin's (now funny) book called ``Materialism and Empirio-Criticism'' where he tries to denounce the works of Mach and Poincare that he could not understand.} For example, celestial mechanics is marvelous because we can describe the observations of planets for several centuries just knowing their masses, positions and velocities at one fixed moment.

The situation becomes more complex if we consider probabilistic models. Imagine that we observe $n$ coin tosses and conclude that they can be described by the Bernoulli model (independent experiments with probability success $1/2$ that lead to the uniform distribution on all $n$-bit strings). This model does not provide a short description of the observed sequence (and, moreover, says that such a description does \emph{not} exist). If so, what is the value of this model? 

The two-part descriptions approach answers this question as follows: 
\begin{itemize}
\item this model is simple; 
\item the two-part description based on this model is close to an optimal one (almost all the time).
\end{itemize}

The second claim, of course, cannot be proven: having written down the results of the coin tossing, we cannot compute the Kolmogorov complexity and check that it is close to the length of the sequence: the Kolmogorov complexity function is not computable. But at least this model can be falsified: if someone comes up with a short description for the results of coin tossing, we are forced to abandon this model. (One may recall Popper here.)

A technical remark: following Kolmogorov, we consider (for simplicity) finite sets as models. In other words, we restricted ourselves to uniform distributions on finite sets as models. But one can easily extend this approach to arbitrary computable probability distributions $P$ (with computable or even rational values), and consider $K(P)+\log_2(1/P(x))$ as the ``length'' of the two-part description provided by model $P$ for some $x$. This ``length'' always exceeds $K(x)$ (up to logarithmic terms), and one may look at how big the excess is. This more general class of models (computable distributions instead of finite sets) does not give essentially new models: one can convert an arbitrary distribution into an uniform one with logarithmic loss of ``quality'', so if we are interested in results with $O(\log n)$ precision for $n$-bit strings, we may consider only finite sets as models.

Several researchers considered two-part descriptions (and corresponding models) at different levels of mathematical precision. Rissanen formulated his ``minimum description length principle'' in his 1978 paper~\cite{1978rissanen}; a more clear exposition (following partially the ideas of Kolmogorov) is provided in his 1999 paper~\cite{1999rissanen}.\endnote{%
Here is the paragraph about two-part description from Rissanen's 1999 paper~\cite{1999rissanen}:
 
       \begin{quote}
First, a `summarizing property' of data may be formalized as a subset $A$ where the data belongs along with other sequences sharing this property. Hence, the property $A$ need not specify the sequence completely. We may now think of programs consisting of two parts, where the first part describes optimally such a set $A$ with the number of bits given by the Kolmogorov complexity $K(A)$, and the second part merely describes $x^n$ in $A$ with about $\log |A|$ bits, $|A|$ denoting the number of elements in $A$. The sequence $x^n$ then gets described in $K(A)+\log|A|$ bits. We may now ask for a set $\hat{A}$ for which $K(\hat{A})$ is minimal subject to the constraint that for an increasing length sequence $x^n$, $K(\hat{A})+\log|\hat{A}|$ agrees with the Kolmogorov complexity $K(x^n)$ to within a constant not depending on $n$. The set $\hat{A}$, or its defining program, may be called Kolmogorov's \emph{minimal sufficient statistic} for the description of $x^n$. The bits describing $\hat{A}$ are then the `interesting' bits in the program (code) for $x^n$ while the rest, about $\log|\hat{A}|$ in number, are non-informative noise bits.
        \end{quote}
A technical remark: the first part should describe $A$ as a finite object (list of all elements); if we consider the program that enumerates $A$ but never terminates explicitly, as a valid description of $A$, then the definition becomes trivial (see the discussion in~\cite{1999shen}).        
}

See also the survey of Grunwald~\cite{2004grunwald} and other papers~\cite{2000vitanyi-li,2000gao-li-vitanyi,2006vitanyi}\endnote{%
Here is the description of ``crude, two-part version of MDL [minimum description length] principle (informally stated)'' taken from~\cite[p.~11]{2004grunwald}.
\begin{quote}
The best point hypothesis $H$ $\langle\ldots\rangle$ to explain the data $D$ is the one which minimizes the sum $L(H) + L(D\cnd H)$, where
\begin{itemize}
\item
$L(H)$ is the length, in bits, of the description of the hypothesis; and
\item
$L(D\cnd H)$ is the length, in bits, of the description of the data when encoded
with the help of the hypothesis.
\end{itemize}

$\langle\ldots\rangle$

We can typically find a very complex point hypothesis (large $L(H)$) with a very good fit (small $L(D\cnd H)$). We can also typically find a very simple point hypothesis (small $L(H)$) with a rather bad fit (large $L(D\cnd H)$). The sum of two description lengths will be minimized at a hypothesis that is quite (but not too) `simple', with a good (but not perfect) fit.
\end{quote}

A technical correction is needed here: the minimal sum is always achiever from the case of singleton description ($H$ allows only $D$); our goal is not just to find the minimal sum, but to find a two-part description of almost the same quality with the simple first part.
}

Some version of the two-description approach was put forward in 1968 (i.e., before Kolmogorov) by Wallace and Boulton~\cite{1968wallaceboulton}\endnote{%
At least one could try to interpret in this way what is said in~\cite[p.~186]{1968wallaceboulton}:
\begin{quote}
\ldots if the observed distribution of the $S$ given points is markedly non-uniform, the economics of Shannon's theorem can be realized in part by adding to the $S$ messages a further message which approximately describes the average density distribution of the $S$ points, and using this approximate non-uniform distribution as the basis for encoding the attribute messages. The message length needed to describe the average density distribution must be considered as adding to the lengths of the messages now used to encode the attributes of the things, because the receiver of the composite message has no \emph{a priori} knowledge of the distribution.
\end{quote}
Here the first part of the description provides an approximate point distribution (for independent identically distributed variables), and the second part corresponds to the minus logarithm of the probability of a point according to this distribution.
}

Now we switch to the third approach suggested by Kolmogorov (under the name of stochasticity), and then we will explain in which sense all three approaches are equivalent (provide the same curve).

\section{Stochasticity}

In 1981 Kolmogorov announced a seminar ``Description complexity and computational complexity'' (it still exists, now in online format, see~\cite{2023kolmogorovseminar}) at the Mathematics Department of Moscow State (Lomonosov) University. (Just before that he became the chairman of the Mathematical Logic Division of this department --- now called Matematical Logic and Algorithms Theory Division.) Kolmogorov opened this seminar and delivered several talks; some notes from these talks (October 28 and November 26, 1981, and October 14, 1982) exist.\endnote{%
Since almost no information about the mathematical interests of Kolmogorov in his last years (in 1981 he was already ill, hardly can walk, his vision deteriorated,so  it was quite difficult for him to give a talk) is available, we reproduce here (with minimal editing) the notes made by one of the authors (A.S.) during these talks.

\medskip
\textbf{October 28, 1981}.

Define
\[
H_B^k(x)=\min\{l(p) \mid B(p)=x \ \text{and computational complexity of $B(p)$ does not exceed $k$}\}
\]
\begin{center}
\includegraphics[scale=1]{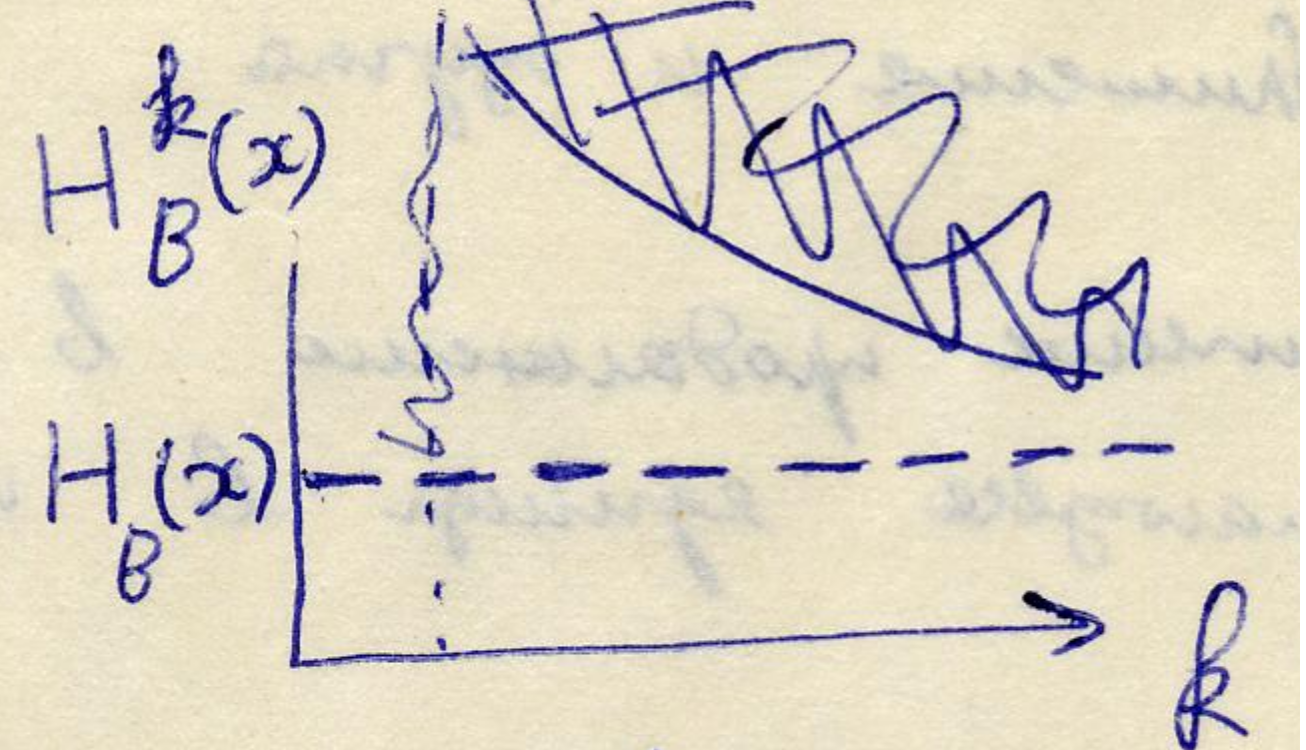}
\end{center}
An example when $H_B^k$ decreases as $k$ increases: old mathematical tables contained some interpolation instructions that make the table shorter (while providing the same information as the original table); the price is that now we need to perform some computations.
\medskip
\[
| I(x\colon y)-I(y\colon x) | = O(\log H(x,y))
\]
We should not be afraid of logarithms (as well as $O(1)$ terms that we have anyway).

Why do we need other versions of algorithmic entropy?

[Other topics for the seminar:]

Boolean complexity of integer and matrix multiplication (Karatsuba, Toom, Strassen).

Assume that we want to compute some function of type $[-1,1]\to[-1,1]$ with precision $\varepsilon$, and are interested in the size of the corresponding Boolean circuit. Let $A_r$ be the functions that have an analytic extension to an ellipse; one can get some lower bound using $\varepsilon$-entropy (Kolmogorov, Tikhomirov). Usually the worst-case complexity for some class is close to the complexity of the majority of elements of this class. But how can we provide an individual example of a hard-to-compute function?

Consider Boolean functions of some fixed complexity; what is the maximal complexity of the inverse function.

\bigskip

\textbf{November 26, 1981}

\medskip

Let $x_1\ldots x_N$ be a bit sequence with $M$ ones, then the complexity is at most $NH(p)$, where $p=M/N$. Assume that for a given sequence it is close to this bound.

A finite object $\omega$ is $(\alpha,\beta)$-stochastic (where $\alpha$ and $\beta$ are some parameters that are small compared to the complexity of $\omega$), if there exists a finite set $\Omega$ such that 
\begin{itemize}
\item $\omega \in  \Omega$;
\item $K(\Omega)\le \alpha$;
\item $K(\omega)\ge K(\Omega)-\beta$ [a typo: there should be  $\log |\Omega|-\beta$ in the right hand side]
\end{itemize}
Let $N(k,\alpha,\beta)$ be the number of $(\alpha,\beta)$-stochastic objects of complexity at most $k$. Typical cases: $\alpha,\beta=\varepsilon k$ or $\alpha,\beta=\sqrt{k}$. What can we say about the asymptotics of $N(k,\alpha,\beta)/2^k$ as $k\to\infty$?  Is it true that $N(k,\alpha,\beta)/2^k \to 0$, if $k$ increases and $\alpha/k, \beta/k \to 0$?

\bigskip
 
\textbf{October 14, 1982}

\medskip

Let $x_1\ldots x_n$ be a bit sequence of length $n$ that contains $m$ ones. The maximal complexity of such a sequence is about $nH(m/n)$ where $H$ is the Shannon entropy function. Let us consider \emph{stochastic} sequences whose complexity is at least this bound minus some $\beta$.

Let us apply some simple selection rule to select $n_1$ elements from this sequence; then $n_2=n-n_1$ elements remain non-selected. Compute the frequencies of ones in both selected and non-selected elements. Under natural assumptions ($n_1$ and $n_2$ are not too small) these frequencies $m_1/n_1$ and $m_2/n_2$ will be close to $m/n$. Indeed, the complexity of the original sequence is bounded by $n_1H(m_1/n_1)+n_2H(m_2/n_2)$ plus the complexity of the selection rule, and at the same time it should be close to $nH(m/n)$. [Then we may use convexity arguments.] What will happen for the dependent selection? How can we reconstruct properties of Markov chains based on maximal complexity assumptions? Sequences with given transition probabilities. Which are stochastic? Prove the known (in probabilistic setting) properties of Markov chains for those sequences.

Let $X_1,\ldots,X_N$ be a finite sequence of integers with fixed mean and variance. Consider a sequence of maximal complexity in this class. Will it have Gaussian distribution? Probably yes: maximizing 
\[
\int_0^\infty f(x)\ln f(x)
\]
with restrictions
\[
\int_0^\infty f(x)\,dx=1, \quad \int_0^\infty x^2 f(x)\, dx =1,
\]
we get a function of type $ce^{-ax^2}$. [See~\cite{1987asarin,1987asarin2}.]

\begin{center}
 $*$ \qquad $*$ \qquad $*$
\end{center}

It looks like these talks were essentially the last mathematical talks given by Kolmogorov. Between them Kolmogorov also gave a talk at the USSR -- Japan symposium on probability theory and mathematical statistics (August 23--29, 1982); the proceedings of this symposium were published in 1983~\cite{1983proceedings}. In his letter (August 6, 1982) Kolmogorov writes to the organizers:

\begin{quote}
Dear colleagues,

I thank you for the invitation to give a talk at the symposium. However, when we had a preliminary discussion about that, I did not have in mind that the working language of the symposium is only English. In general this looks reasonable and I don't want the organizers to make an exception for me. So I have to decline your invitation. If you still want to honor my position as a doyen of probability theory, one of the options is to let me deliver an introductory speech (15 minutes in Russian with 15-minutes for the English translation); I will then say something meaningful about the tendencies in the development of our subject. Anyway, I will try to be useful for the symposium. I will arrive on August 22 by train.
\end{quote}

Kolmogorov's introductory speech was recorded (now the recording seems to be lost) and then transcribed for publication~\cite{1983kolmogorovconf} (this was done by A.K.~Zvonkin and A.~Shen as requested by A.A.~Novikov). However, the recording was very bad, the voice of Kolmogorov (and even the translator) was almost inaudible, so the published text is more a reconstruction of what Kolmogorov could say than the transcript --- the reconstruction based on the seminar talks of Kolmogorov and his earlier publications. Here is the fragment about the notion of a random (typical) element of a finite set (here $A$ is the optimal description method fixed when we define conditional complexity):
\begin{quote}
Now we can define the concept of a ``random'', or, to be more precise, $\Delta$-random object in a given finite set $M$ (here $\Delta$ is a number). Namely, we shall say that $X\in M$ [a typo: $\in$ sign is missing in the text] is $\Delta$-random in $M$ if
\[
K_A(X\cnd M) \ge \log_2 |M| - \Delta,
\]
[another typo in the text: $Y$ is printed instead of $M$] where $|M|$ [one more typo in the text: $M$ is printed instead of $|M|$, the cardinality of $M$] denotes the number of elements in $M$. We shall call \emph{random in $M$} the $\Delta$-random objects in $M$, $\Delta$ being comparatively small. Thus we receive the definition of a random finite object which can be regarded as a final one.
\end{quote}

When in 1986 the first world Bernoulli Society congress was organized (the opening day was September 8, 1986). Kolmogorov was too ill to come, and V.A.~Uspensky read his welcome remarks:

\begin{quote}
Dear ladies and gentlemen! 

Allow me to welcome you today to the opening of the Congress.

It is significant to me that the Society that has taken the name Bernoulli, a Society uniting specialists in just one field of mathematics --- probability theory and mathematical statistics --- has succeeded in organizing a conference of its fellow members so representative that it is comparable to international mathematics congresses. But if one thinks about it, then one can find an explanation for this seemingly paradoxical phenomenon.

James [Jacob] Bernoulli, one of the eminent members of the Bernoulli family, has entered the pages of the history of science by virtue of his many achievements. But two of his credits should be mentioned especially. He is the father of the science of probability theory having obtained the first serious result known everywhere as Bernoulli's theorem. But apart from this, it should not be forgotten that he was essentially also the father of combinatorial analysis. He used the elements of this discipline to prove his theorem but he delved into the field of combinatorial analysis considerably further discovering in particular the remarkable sequence of numbers which now bear his name. These numbers are encountered continually in scientific investigations right down to our time.

We all feel that one of the basic requirements of mathematics that is evident at present is the investigation of very complex systems. And this complexity on the one hand is very closely related to randomness and on the other --- it necessitates in some measure an extension of combinatorial analysis itself. All this gives hope that as time passes the Bernoulli Society will increase its influence more and more in the mathematical world. I wish the participants of the Congress all of the very best.

\end{quote}

After this welcome remarks, there was a plenary session of the Congress (in an opera theater) where Uspensky (speaking from the orchestra pit and not seeing the screen with his presentation) gave a talk that was declared as a joint talk of Kolmogorov and Uspensky (a paper based on this talk that appeared later~\cite{1987kolmogorov} had the same authors). However, one should have in mind that (as it is written in the Introduction), ``while the main contents of this paper  (especially of its two first chapter) is based on the ideas and publications of A.N.~Kolmogorov, the first author [Kolmogorov], who is a teacher of the second author [Uspensky] was not able to look at the final version of the text''. So the definition of stochasticity given in this paper reflects the discussions on the seminar after Kolmogorov's talks.
}

In the first of these talks Kolmogorov reproduced the definition of time-bounded complexity and illustrated it by interpolation rules that quite often appeared in old mathematical tables: they decrease the size of the table but make its use more time-consuming.

In the second talk he suggested the definition of a stochastic object (an object that can be explained statistically). Later it was reproduced in~\cite[p.~1337]{1983shen} (with some lower and upper bounds for the fraction of non-stochastic objects of given size):
\begin{quote}
\textbf{Definition} (A.N.~Kolmogorov). Let $\alpha$ and $\beta$ be some natural numbers. A number $x$ will be called $(\alpha,\beta)$-\emph{stochastic} if there exists a finite set $A\subset \mathbb{N}$ such that 
\[
x\in A, \qquad K(A)\le \alpha, \quad K(x) \ge \log_2\#A-\beta,
\]
where $\#A$ denotes the number of elements in $A$,
\end{quote}
and $K$ stands for Kolmogorov complexity.

The idea of this definition becomes clear if we compare it with Kolmogorov's remark from~\cite{1965kolmogorov}:

\begin{quote}
Another topic beyond the scope of this note is the application of constructions developed in \S3 to a new basis for probability theory. Roughly speaking, this application can be described as follows. Consider a finite $M$ that has a very large number $N$ of elements and can be described by a program whose length is negligible compared to $\log_2 N$. Then almost all elements $x$ of $M$ have complexity $K(x)$ close to $\log_2 N$. The elements of $M$ that have this complexity can be considered as ``random'' elements of $M$.
\end{quote}

Then, in the definition of stochasticity, given some $x$, we look for a simple set $A$ (condition $K(A)\le\alpha$) that contains $x$ such that $x$ is a ``random'' element of $A$ (condition $K(x)\ge \log_2 |A|-\beta$).

If the complexity of $A$ is not negligible, it is natural to add $A$ to the condition in the second part of this definition, as it is done in~\cite{1987kolmogorov}:
\begin{quote}
The \emph{randomness deficiency} of some element $y$ in the set $M$ is defined as 
\[
d(y\cnd M)=\log_2 \#M- H(y\cnd M),
\]
where $\#M$ stands for the cardinality of the set $M$. $\langle\ldots\rangle$

We say that $y\in M$ is \emph{$\Delta$-random} in $M$ if $d(y\cnd M)\le \Delta$. Then ``sufficiently random elements of $M$'' can be defined as $\Delta$-random elements of $M$ for sufficiently small  $\Delta$. This definition was suggested in~[35]  $\langle\ldots\rangle$

A natural question arises: are there ``absolutely non-random objects'', i.e., objects that have a large randomness deficiency in all simple sets that contain them? The answer is positive: such objects do exist. $\langle\ldots\rangle$

This definition has the following motivation. A statistician gets some experimental data $y$ and wants to explain this data. This means to provide some statistical model (hypothesis), in other words, to find some set of ``possible outcomes'' of the experiments that contains $y$ such that $y$ looks like a typical (random) element. In mathematical terms, we look for a simple finite set $A$ containing $y$ such that $y$ has a small randomness deficiency in $A$.
\end{quote}

Here [35] is the Kolmogorov's 1982 talk~\cite{1983kolmogorovconf}.

It seems that no participants of the 1981--1982 Kolmogorov seminar were aware of his earlier talks where he introduced the structure function (or have seen his 1974 abstract~\cite{1974umn}). So nobody asked Kolmogorov what he thinks about the connection of his new definition (that did not seem to appear earlier in his talks and publications) and the structure functions. Soon after that some papers written by the seminar participants appeared~\cite{1983shen,1985vyugin,1987vyugin,1999vyugin} that studied the notion of $(\alpha,\beta)$-stochasticity without mentioning the structure function in any way.

In a strange way some participants of the Tallinn 1973 conference (where the structure function was defined) relate Kolmogorov's talk there with the notion of stochasticity. In~\cite[page 858]{1989CGG} we find the following remark: 
\begin{quote}
Kolmogorov introduced (see [K462]) the quantity 
\[
d(x\cnd S) = \log |S| - K(x\cnd S)
\]
as the defect of randomness of a string $x$ with respect to set $S$. $\langle\ldots\rangle$

At a Tallin[n] conference in 1973, Kolmogorov proposed a variant of the function
\[
\delta_x(k) = \min_S\{d(x\cnd S)\colon K(S)<k, x\in S\},
\]
considering it an interesting characteristic of the object~$x$. If $x$ is some sequence of experimental results, then the set $S$ can be considered to be the extraction of those features in $x$ that point to nonrandom regularities. At the point $k^*(x)$ where the decreasing function $\delta_x(k)$ becomes zero (or less than some constant agreed on in advance), we can say that it is useless to explain $x$ in greater detail than by giving the set $S^*$ such that $d(x\cnd S^*)=\delta_x(k^*)$. Indeed, the added explanation would be as large as the number of extra bits it accounts for. The set $S^*$ expresses all the relevant structure in the sequence $x$, the remaining details of $x$ being conditionally maximally random. For example, $S^*$ would describe the Mona Lisa up to brush strokes, and $k^*$, the length of description of $S^*$, is the ``structural complexity of $x$''.

The set $S^*$ plays the role of a \emph{universal minimal sufficient statistic} for $x$. 

\end{quote}

Here [K462] refers to~\cite{1983kolmogorovconf} (though the randomness deficiency is not explicitly defined there). Maybe the words ``a variant of'' mean that the authors of the survey~\cite{1989CGG} prefer to use notions that were introduced later by Kolmogorov? 

Anyway, a natural question arises. Kolmogorov suggested two criteria saying when a finite set of strings $S$ is a good model (statistical hypothesis) for a string $x$:
\begin{itemize}
\item if a two-part description of $x$ based on $S$ is close to the optimal description of $x$, i.e., the ``optimality deficienty''
\[
\delta(x\cnd S) = K(S)+\log_2 \#S - K(x);
\]
is small
\item if the randomness deficiency
\[
d(x\cnd S) = \log_2 \#S - K(x\cnd S).
\]
is small.
\end{itemize}
Are these two conditions equivalent? It is easy to see (using the formula for the complexity of a pair) that the optimality deficiency (with logarithmic precision) exceeds the randomness deficiency by $K(S\cnd x)$:
\[
\delta(x\cnd S) = d(x\cnd S)+K(S\cnd x)+O(\log n),
\]
(for string $x$ and set $S$ of the complexity at most $n$). This difference is negligible in the case (that is most interesting from the philosophical viewpoint) when $S$ is simple. Then its conditional complexity given $x$ (or without this condition) is small. But in general case these two conditions are not equivalent, as noted in~\cite{2001gacs-tromp-vitanyi}. Indeed, consider two independently chosen random $n$-bit strings $x$ and $y$ and consider the set $S=\{x,y\}$ as a model for one of them (say, $x$). Then $\delta(x\cnd S)$ is about $n$ (the complexity of $S$ is close to $2n$, and this is the main part of a two-part description based on $S$) while $d(x\cnd S)$ is close to zero (since $\log_2\#S$ is only $1$).

A more interesting question was posed in~\cite{2001gacs-tromp-vitanyi} (as an open question). Assume that for some $x$ we are able to find a model (set $S$) with some values of $K(S)$ and $d(x\cnd S)$.  Is it possible to find (another) set $S'$ such that $K(S')$ does not exceed much the complexity $K(S)$ and at the same time $\delta(x\cnd S')$ does not exceed much $d(x\cnd S)$? This question was soon solved in~\cite{2004vereshchaginvitanyi} in a positive direction. So it became known that two approaches suggested by Kolmogorov (the structure function and $(\alpha,\beta)$-stochasticity) lead to the same curve (with logarithmic precision). Moreover, it turned out that the same curve can be defined in terms of resource-bounded complexity. In the next section we formulate these results more explicitly and describe their history.

\section{Three equivalent approaches}

As we said, all three approaches suggested by Kolmogorov are equivalent. Each of them, starting from a string $x$ of length $n$, gives some plane curve, and these three curves turn out to be the same for the appropriately chosen coordinates (with $O(\log n)$ precision).

Let us repeat the definitions of these three curves. Let $x$ be an $n$-bit string.
\begin{itemize}
\item\textbf{Resource-bounded complexity}.

Let $B(k)$ be the maximal natural number of complexity at most $k$. Consider the function
\[
k \mapsto K^{B(k)}(x) - K(x)
\]
(showing how much the resource-bounded complexity exceeds the unbounded one). The same curve (in a different coordinate system) appears in the definition of logical/computational depth.

\item\textbf{Structure function}
Consider the function 
\[
k \mapsto k + H_k(x) - K(x), 
\]
where $H_k(x)$ is the Kolmogorov structure function (also denoted by $\Phi_k$, see Section~\ref{sec:structure}). In other words, for every $k$ we consider the minimum value of the optimality deficiency $\delta(x\cnd S)$ of all finite sets $S$ of complexity at most $k$. The same curve in different coordinates appears in the definition of sophistication.

\item \textbf{Stochasticity}
Consider the function 
\[
k \mapsto \min \{l \colon \text{$x$ is $(k,l)$-stochastic}\}.
\]
In other words, for every $k$ we consider the minimal randomness deficiency $d(x\cnd S)$ over all finite sets $S$ of complexity at most $k$.
\end{itemize}

\smallskip
\begin{center}
\fbox{\fbox{\parbox{0.48\textwidth}{%
\textbf{Main theorem of algorithmic statistics}.\\ For every $n$-bit string $x$ these three functions coincide with $O(\log n)$ precision.}}}
\end{center}
\smallskip
When speaking about $O(\log n)$-precision, we allow $O(\log n)$-changes in both coordinates. More formally, we consider three sets that are upper-graphs of these three functions (the second coordinate is at least the value of the function on the first coordinate) and claim that each of them is covered by the $O(\log n)$-neighborhood of every other one.\endnote{%
There is one more equivalent description of this curve --- not related to Kolmogorov --- that goes back to~\cite{2001gacs-tromp-vitanyi}. Assume that $m\ge K(x)$; consider an algorithm that, given $m$, enumerates all strings of complexity at most $m$. The string $x$ appears in this enumeration. Look how many strings appear \emph{after} $x$; assume that there are about $2^s$ of them. Then consider the point $(m-s,s)$ on the curve (that corresponds to a two-part description of complexity $m-s$ and size $s$).
}

As we have said already, the equivalence between the second and third approaches was established in~\cite{2004vereshchaginvitanyi}. As to the first approach, the equivalence is proven in~\cite{2017abst}. This paper seems to be the first one where this equivalence result is established in that form (the coincidence of two curves with logarithmic precision), though some related results were obtained earlier; for example, already in~\cite{1987koppel} Koppel claimed that some related numerical quantities (but not curves) defined for infinite sequences are close to each other\endnote{%
In Kolmogorov's opinion, we should first of all consider finite objects, as these considerations could lay a theoretical foundation for practical use of probability theory. The infinite objects, from his viewpoint, were kind of decorations that would be nice but not philosophically important. He wrote to Uspensky (June 2, 1983):
\begin{quote}
Dear Vladimir Andreevich!

Of course I cannot have any objections if you want to write a paper with Shen and Semenov about different algorithmic definitions of a notion of an infinite random sequence. However, as to the publication of this paper in the \emph{Uspekhi matematicheskikh nauk} [Russian Mathematical Surveys] journal, I have some reservations. I do not consider the problem of defining infinite random sequences as a central problem. To understand how I see the general picture, please look at

1. \S2 of my book ``Foundations on the theory of probability'';

2. my paper in the collection ``Mathematics, its methods and role'';

3. my paper in \emph{Sankhya} ``On the tables of random numbers'' [see~\cite{1963kolmogorov}];

4. the text of my Nice talk (attached) [see~\cite{1970kolmogorov}].

I speak about infinite sequences in section 6 and 8 of talk in Nice. The infinite random sequences remain a nice appendix that is not really necessary. This is also the case for the traditional exposition of probability theory. There are two types of limit theorems: for infinite and finite series of trials, and the second one is (rightly) considered as being more close to the practical needs. 
\end{quote}

Since then the randomness theory for infinite sequences developed a lot, mostly due to recursion theorists who came to this field and found many surprising connections between general recursion theory and algorithmic information theory (e.g., a very interesting notion of K-trivial Turing degrees). A lot of nice (and often difficult) results of this type can be found in~\cite{2010downeyhirschfeldt,2012nies}. However, they are interesting mostly from the mathematical point of view and rarely give something to the philosophical foundations of probability theory. (Still one should note that Schnorr's ideas about martingale randomness and algorithmic information theory in general were one of the sources of inspiration of the game approach to probability theory~\cite{2001vovkshafer}.)
}
; see also~\cite{2007antunesfortnow,2010bauwens}

What Kolmogorov thought about these equivalences is not clear. Did he know (at least as a conjecture) that all three definitions will turn out to be equivalent? Probably we will never know this. But the fact that he has suggested, on three different occasions, three ways to measure the ``complexity behavior'' of a string not only by some number, but by some curve, that later turned out to be equivalent, is in itself a remarkable achievement that can be rightfully called ``the last discovery of Kolmogorov''.

\begin{center}
$*$\qquad $*$ \qquad $*$
\end{center}

In this paper we tried to follow the development of Kolmogorov's idea, and further results in the field of algorithmic statistics are beyond the scope of our survey. The proofs of the equivalence between three approaches, as well as some subsequent results, can be found in~\cite{2016vereshchaginshen}; they are also discussed (without proofs) in a shorter survey~\cite{2015vereshchaginshen}. See also the last chapter of~\cite{2013SUV}.

The authors are deeply grateful to Albert Nikolaevich Shiryaev, in particular, for his constant effort to keep the memories about Andrei Nikolaevich Kolmogorov, and to all participants of seminars at Moscow University Mathematics Department (Logic and Algorithms Division), LIF (Marseille) and LIRMM (Montpellier) where the notions mentioned in this survey were discussed.

Authors acknowledge with gratitude the influence of Vladimir Andreevich Uspensky (1930--2018) and Andrei Albertovich Muchnik (1958--2007) who will remain in the authors' memory.

\printendnotes[custom]


\begin{thebibliography}{9}
\raggedright

\bibitem{1960solomonoff}
R.J.~Solomonoff, \emph{V-131. A preliminary report on a general theory of inductive inference}. February 4, 1960. Zator company, $140\frac{1}{2}$ Mount Auburn Street, Cambridge 28, Mass. Contract AF 49(638)-376. Air Force Office of Scientific Research, Air Research and Development Command, United States Air Force, Washington 25, D.C., available at~\url{http://raysolomonoff.com/publications/rayfeb60.pdf}

\bibitem{1963kolmogorov}
A.N.~Kolmogorov, On tables of random numbers, \emph{Sankhya. The Indian Journal of Statistic}, 1963, Series A, \textbf{25}(4), 369--376. \url{https://www.jstor.org/stable/25049284}. Reprinted in \emph{Theoretical Computer Science}, \textbf{207}, 387--395 (1998), \url{https://www.sciencedirect.com/science/article/pii/S0304397598000759}

\bibitem{1964-1solomonoff}
R.J.~Solomonoff, A Formal Theory of Inductive Inference. Part I, \emph{Information and control}, \textbf{7}(1), 1--22 (1964), \url{https://doi.org/10.1016/S0019-9958(64)90223-2}

\bibitem{1964-2solomonoff}
R.J.~Solomonoff, A Formal Theory of Inductive Inference. Part II, \emph{Information and control}, \textbf{7}(2), 224--254 (1964), \url{https://doi.org/10.1016/S0019-9958(64)90131-7}

\bibitem{1965kolmogorovtalk}
\rus{А.\,Н.\,Колмогоров. Стенограмма доклада <<Понятие ``информация'' и основы теории вероятностей>>. В книге: \emph{Колмогоров и кибернетика} под редакцией Д.\,А.\,Поспелова, Я.\,И.\,Фета. Новосибирск, ИВМиМГ СО РАН, 2001. -- 159 с. (Вопросы истории информатики. Выпуск 2), с.~118--142,} \url{https://archive.org/details/kolmogorov-1965talk}

\bibitem{1965kolmogorov}
\rus{А.\,Н.\,Колмогоров, Три подхода к определению понятия <<количество информации>>, \emph{Проблемы передачи информации}, том 1, вып.~1, 3--7 (1965), \url{https://www.mathnet.ru/rus/ppi68}. (Поступила в редакцию 9 января 1965 года)}
English version: Three approaches to the quantitative definition of information, International Journal  Comput. Math., \textbf{2}, 157--168.

\bibitem{1968umn}
\rus{В Московском математическом обществе. Заседания Московского математического общества. \emph{Успехи математических наук}, XXIII, вып. 2 (140), 1968, март--апрель, с.~201}, \url{https://www.mathnet.ru/rus/rm5615}

\bibitem{1968wallaceboulton}
S.C.~Wallace, D.M.~Boulton, An information measure for classification, \emph{Computer Journal}, \textbf{11}(2), 185--194, \url{https://academic.oup.com/comjnl/article/11/2/185/378628}.

\bibitem{1969kolmogorov}
\rus{А.\,Н.\,Колмогоров, К логическим основам теории информации и теории вероятностей, \emph{Проблемы передачи информации}, том 5, вып. 3, 3--7 (1969), \url{https://www.mathnet.ru/rus/ppi1805}.} English version: Andrei N.~Kolmogorov, Logical Basis for Information Theory and Probability Theory, IEEE Transactions on information theory, IT-14,  no.~5, 662--664 (September 1968, received December 13, 1967, based on an invited lecture at the International Symposium on Information Theory, San Remo, Italy, September, 1967.)

\bibitem{1970kolmogorov}
\rus{А.\,Н.\,Колмогоров. Комбинаторные основания теории информации и исчисления вероятностей. \emph{Успехи математических наук}, том 38, выпуск 4(232), июль-август 1983, 27--36, \url{https://www.mathnet.ru/rus/rm2940}}. Original Kolmogorov's manuscript of 1970: \url{https://archive.org/details/kolmogorov83-manuscript}

\bibitem{1972umn} 
\rus{В Московском математическом обществе. Заседания Московского математического общества, \emph{Успехи математических наук}, XXVII, выпуск 2, с.~159}, \url{https://www.mathnet.ru/rus/rm5033}

\bibitem{1973fine}
Terrence L.Fine, Photo of Kolmogorov giving a talk in Tallinn, \url{https://commons.wikimedia.org/wiki/File:Kolm_complexity_lect.jpg}. Description: ``taken at the 1973 Soviet Information Theory Symposium (may not be the exact title) held in Tallinn, Estonian SSR. Kolmogorov delivers his talk. A.M. Yaglom is also on the picture.''

\bibitem{1973dynkin}
Audio recording of some talk of Kolmogorov from Evgeny B.~Dynkin collection (\url{https://ecommons.cornell.edu/handle/1813/17350}). In this collection this recording is annotated as Kolmogorov's interview given to Dynkin in 1975; however, this is not an interview but a talk at some international conference (one could hear translator's voice in background). Maybe this is Kolmogorov's talk at Tallinn's conference in 1973. The version with Russian subtitles: \url{https://www.youtube.com/watch?v=GN519NP2JXI}

\bibitem{1974umn}
\rus{В Московском математическом обществе. Заседания Московского математического общества. \emph{Успехи математических наук}, XXIX, выпуск 4 (178), 1974, июль--август.} \url{https://www.mathnet.ru/rus/rm7215}

\bibitem{1977chaitin}
Gregory Chaitin, Algorithmic information theory, \emph{IBM Journal of Research and Development}, \textbf{21}(4), 350--359 (1977), \url{https://dl.acm.org/doi/10.1147/rd.214.0350}

\bibitem{1978rissanen}
J.~Rissanen, Modeling by shortest data description, \emph{Automatica}, \textbf{14}(5), 465--471, \url{https://www.sciencedirect.com/science/article/abs/pii/0005109878900055}, see also \url{https://msol.people.uic.edu/ECE531/papers/Modeling%20By%20Shortest%20Data%20Description.pdf}.

\bibitem{1979adleman}
Leonard M.~Adleman, \emph{Time, Space and Randomness}, MIT report LCS/TM-131, April, 1979. Available at \url{https://people.cs.rutgers.edu/~allender/papers/adleman.time.space.randomness.pdf}

\bibitem{1983proceedings}
\emph{Probability Theory and Mathematical Statistics. Proceedings of the Fourth USSR--Japan Symposium, held in Tbilisi, USSR, August 23--29, 1982}, edited by Jurii V.~Prokhorov, Kiyosi Itô (Lecture Notes in Mathematics, v.~1021). Springer, 1983, \url{https://doi.org/10.1007/BFb0072896}. 

\bibitem{1983kolmogorovconf}
A.N.~Kolmogorov, On logical foundations of probability theory.  In~\cite{1983proceedings}, 1--5. (1983)

\bibitem{1983gacs}
P\'eter G\'acs, On the relation between descriptional complexity and algorithmic probability, \emph{Theoretical Computer Science}, \textbf{22}(1--2), 71--93, see \url{https://doi.org/10.1016/0304-3975(83)90139-1}.

\bibitem{1983sipser}
Michael Sipser, A complexity theoretic approach to randomness,  \emph{STOC'83: Proceedings of the 15th annual ACM symposium on Theory of computing}, 330--335 (December 1983), see \url{https://dl.acm.org/doi/10.1145/800061.808762}.

\bibitem{1983shen}
\rus{А.\,Х.\,Шень, Понятие $(\alpha,\beta)$-стохастичности по Колмогорову и его свойства. \emph{Доклады Академии наук СССР}, \textbf{271}(6), 1337--1340 (1983),} \url{https://www.mathnet.ru/rus/dan9984}.

\bibitem{1985cover}
Thomas M.~Cover, Kolmogorov Complexity, Data Compression, and Inference. In: \emph{The Impact of Processing Techniques on Communications}, edited by J.K.~Skwirzhynski, Martinus Nijhoff Publishers, 1985, \url{https://isl.stanford.edu/~cover/papers/paper65.pdf}. 

\bibitem{1985vyugin}
\rus{В.\,В.\,Вьюгин, О нестохастических объектах, \emph{Проблемы передачи информации}, \textbf{21}, выпуск 2, 3--9 (1985)}, \url{https://www.mathnet.ru/rus/ppi979}.

\bibitem{1986longpre}
Luc Longpr\'e, \emph{Resource Bounded Kolmogorov Complexity, A Link Between Computational Complexity and Information Theory}, Ph.D thesis, TR 86-776, August 1986, Department of Computer Science, Cornell University, Ithaca, NY14853, see \url{https://ecommons.cornell.edu/handle/1813/6616}.

\bibitem{1987kolmogorovspeech}
\rus{А.\,Н.\,Колмогоров. Приветствие участникам Первого Всемирного конгресса Общества Бернулли (зачитанное В.\,А.\,Успенским). \emph{Теория вероятностей и её применения},} \textbf{32}(2), 218, \url{https://www.mathnet.ru/rus/tvp1408}.

\bibitem{1987kolmogorov}
\rus{А.\,Н.\,Колмогоров, В.\,А.\,Успенский, Алгоритмы и случайность, \emph{Теория вероятностей и её применения}, том 32, выпуск 3 (июль--сентябрь 1987),} 425--455, \url{https://www.mathnet.ru/rus/tvp1437}.

\bibitem{1987asarin}
\rus{Е.\,А.\,Асарин, О некоторых свойствах случайных в алгоритмическом смысле конечных объектов. \emph{Доклады АН СССР}, 1987, том 295, номер 4, 782--785.} \url{https://www.mathnet.ru/rus/dan8027}.

\bibitem{1987asarin2}
\rus{Е.\,А.\,Асарин, О некоторых свойствах $\Delta$-случайных по Колмогорову конечных последовательностей. \emph{Теория вероятностей и её применения}, 1987, том 32, выпуск 3, 556--558.} \url{https://www.mathnet.ru/rus/tvp1450}.

\bibitem{1987vyugin}
\rus{В.\,В.\,Вьюгин, О дефекте случайности конечного объекта относительно мер с заданными границами их сложности, \emph{Теория вероятностей и её применения}}, \textbf{32}(3), 558--563 (1987), \url{https://www.mathnet.ru/rus/tvp1451}. English translation: V.V.~V'yugin, On the defect of randomness of a finite object with respect to measures with given complexity bounds, \emph{Theory Prob. Appl.}, \textbf{32}, 508--512 (1987), DOI:10.1137/1132071 (published in 1988).

\bibitem{1987koppel}
Moshe Koppel, Complexity, Depth and Sophistication, \emph{Complex systems}, \textbf{1}, 1087--1091 (1987), \url{https://www.complex-systems.com/abstracts/v01_i06_a04/}.

\bibitem{1988bennett}
Charles H.~Bennett, Logical Depth and Physical Complexity, in \emph{The Universal Turing Machine --- a Half-Century Survey}, Rolf Herken, editor, Oxford University Press, 1988, 227--257.

\bibitem{1988koppel} 
Moshe Koppel, Structure, in \emph{The Universal Turing Machine --- a Half-Century Survey}, Rolf Herken, editor, Oxford University Press, 1988, 435--452.

\bibitem{1989CGG}
Thomas M.~Cover, Peter Gacs, Robert M.~Gray, Kolmogorov contributions to information theory and algorithmic complexity, \emph{The Annals of Probability}, 1989, vol.~17, no.~3, 840--865. 

\bibitem{1991koppelatlan}
Moshe Koppel, Henri Atlan, 
An Almost Machine-Independent Theory of Program-Length Complexity, Sophistication, and Induction,
\emph{Information Sciences}, \textbf{56}, 23--33 (1991).

\bibitem{1991CT}
Thomas M.~Cover, Joy A.~Thomas, \emph{Elements of information theory}, New York: Wiley, 1991. Second edition was published in 2006.

\bibitem{1992li-vitanyi}
Ming Li, Paul M.B.~Vit\'anyi, Inductive Reasoning and Kolmogorov Complexity, \emph{Journal of Computer and System Sciences}, \textbf{44}(2), 343--384 (1992), \url{https://doi.org/10.1016/0022-0000(92)90026-F} (page 353 is missing). Preliminary version: M. Li and P. M. B. Vitanyi, \emph{Proceedings. Structure in Complexity Theory Fourth Annual Conference, Eugene, OR, USA, 1989},165--185, \url{https://doi.org/10.1109/SCT.1989.41823}

\bibitem{1993longpre}
Luc Longpr\'e, Sarah Mocas, Symmetry of information and one-way functions, \emph{Information Processing Letters}, \textbf{46}(2), 95--100 (1993), \url{https://doi.org/10.1016/0020-0190(93)90204-M}.

\bibitem{1999rissanen}
Rissanen, J.~Hypotheses selection and testing by the MDL principle,
\emph{Computer Journal}, \textbf{42}(4), 260--269 (1999), \url{https://ieeexplore.ieee.org/document/8138703}.

\bibitem{1999shen}
Alexander Shen, Discussion on Kolmogorov complexity and statistical analysis, \emph{Computer Journal}, \textbf{42}(4), 340--342 (1999), \url{https://doi.org/10.1093/comjnl/42.4.340}

\bibitem{1999vyugin}
V.V.~Vyugin, Algorithmic complexity and stochastic properties of finite binary sequences, \emph{Computer Journal}, \textbf{42}(4), 294--317 (1999), \url{https://doi.org/10.1093/comjnl/42.4.294}

\bibitem{2000gao-li-vitanyi}
Qiong Gao, Ming Li, Paul Vit\'anyi, Applying MDL to learn best model granularity,
\emph{Artificial Intelligence}, \textbf{121}, 1--29 (2000), \url{https://homepages.cwi.nl/~paulv/papers/ai00.pdf}.

\bibitem{2001gacs-tromp-vitanyi}
P\'eter G\'acs, John T.~Tromp, Paul M.B.~Vit\'anyi, Algorithmic Statistics, \emph{IEEE Transactions on Information Theory}, \textbf{47}(6), 2443--2463 (September  2001), \url{https://ieeexplore.ieee.org/abstract/document/945257}.

\bibitem{2000vitanyi-li}
Paul M.B.~Vit\'anyi, Ming Li, Minimum Description Length Induction,
Bayesianism, and Kolmogorov Complexity, \emph{IEEE Transactions on Information Theory}, \textbf{46}(2), 446--464 (March 2000), \url{https://homepages.cwi.nl/~paulv/papers/mdlindbayeskolmcompl.pdf}.

\bibitem{2001vovkshafer}
Glenn Shafer, Vladimir Vovk, \emph{Probability and Finance: It's Only a Game}, Wiley, 2001, see also~\url{http://www.probabilityandfinance.com/2001_book/index.html}.

\bibitem{2004grunwald}
Peter Grunwald, \emph{A tutorial introduction to the minimum description length principle}, \url{https://arxiv.org/abs/math/0406077}.

\bibitem{2004vereshchaginvitanyi}
Nikolai Vereshchagin, Paul Vit\'{a}nyi, Kolmogorov’s Structure Functions and Model Selection, \emph{IEEE Transactions on Information Theory}, \textbf{50}(12), 3265--3290 (December 2004). Previous version: 47th FOCS (2002). See also: \url{https://arxiv.org/pdf/cs/0204037.pdf}.

\bibitem{2006AFMV}
Luis Antunes, Lance Fortnow, Dieter van Melkebeek, N.V. Vinodchandran,
Computational depth: Concept and applications, \emph{Theoretical Computer Science}, \textbf{354}(3), 391--404 (April 2006), see \url{https://www.sciencedirect.com/science/article/pii/S0304397505008790}.

\bibitem{2006vitanyi}
Paul M.B.~Vit\'anyi, Meaningful information, \emph{IEEE TRansactions on Information Theory}, \textbf{52}(10), 4627--4626 (October 2006). See also: \url{https://arxiv.org/pdf/cs/0111053v3.pdf}.

\bibitem{2007antunesfortnow}
Luis Antunes, Lance Fortnow, Sophistication revisited, \emph{Theory of Computing Systems}, \textbf{45}, 150--161, \url{https://doi.org/10.1007/s00224-007-9095-5} (2007) Preliminary version: 30th ICALP (2003).

\bibitem{2010bauwens}
Bruno Bauwens, Computability in statistical hypotheses testing, and characterizations of independence and directed influences in time series using Kolmogorov complexity, Ph.D. thesis, University of Ghent, 2010. ISBN 978-90-8578-356-5, available at~\url{https://biblio.ugent.be/publication/1107852}.

\bibitem{2010downeyhirschfeldt}
Rodney G.~Downey, Denis R.~Hirschfeldt, \emph{Algorithmic Randomness and Complexity}, Springer, 2010, \url{https://doi.org/10.1007/978-0-387-68441-3}.

\bibitem{2012nies}
Andre Nies, \emph{Computability and Randomness}, Oxford University Press, 2009, \url{https://global.oup.com/academic/product/computability-and-randomness-9780199652600}.

\bibitem{2013SUV}
\rus{Н.\,К.\,Верещагин, В.\,А.\,Успенский, А.\,Шень, \emph{Колмогоровская сложность и алгоритмическая случайность}. М., МЦНМО, 2013. См. также~\url{https://hal-lirmm.ccsd.cnrs.fr/lirmm-00786255/document}.} English version: A.~Shen, V.A.~Uspensky, N.~Vereshchagin, \emph{Kolmogorov Complexity and Algorithmic Randomness}, American Mathematical Society, 2017, see~\url{https://www.lirmm.fr/~ashen/kolmbook-eng-scan.pdf}, \url{https://hal-lirmm.ccsd.cnrs.fr/lirmm-01803620v1/document}.

\bibitem{2015vereshchaginshen}
Nikolay Vereshchagin, Alexander Shen, 
Algorithmic Statistics Revisited. In: Vovk, V., Papadopoulos, H., Gammerman, A. (eds) \emph{Measures of Complexity. Festschrift for Alexey Chervonenkis}.  Springer, Cham. \url{https://doi.org/10.1007/978-3-319-21852-6_17}, 235--252 (2015), see also~\url{https://arxiv.org/pdf/1504.04950.pdf}.

\bibitem{2016vereshchaginshen}
Nikolay Vereshchagin, Alexander Shen, 
Algorithmic Statistics: Forty Years Later. In: Adam Day, Michael Fellows, Noam Greenberg, Bakhadyr Khoussainov, Alexander Melnikov, Frances Rosamond, editors, \emph{Computability and Complexity. Essays Dedicated to Rodney G. Downey on the Occasion of His 60th Birthday}.  Springer, Cham. \url{https://doi.org/10.1007/978-3-319-21852-6_17} (2016), 669--737,  see also~\url{https://arxiv.org/abs/1607.08077}.

\bibitem{2017abst}
Lu\'\i s Antunes, Bruno Bauwens, Andr\'e Souto, Andreia Teixeira, Sophistication vs Logical Depth, \emph{Theory of Computing Systems}, \text{60}, 280--298 (2017), \url{https://link.springer.com/article/10.1007/s00224-016-9672-6}, see also~\url{https://arxiv.org/pdf/1304.8046.pdf}.


\bibitem{2023kolmogorovseminar}
Some ``Kolmogorov's seminar'' materials (it became an international on-line seminar): \url{https://www.youtube.com/channel/UC20xHyxD6FqItj2N6y3eSZg}.

\end{thebibliography}
\end{document}